\documentclass[12pt]{article}

\usepackage{amsmath,amssymb,amsthm,amscd}
\usepackage[colorlinks, linkcolor=blue, anchorcolor=gree, citecolor=red]{hyperref}
\usepackage[numbers,sort&compress]{natbib}
\usepackage{graphicx}

\makeatletter

\newcommand{\Rmnum}[1]{\expandafter\@slowromancap\romannumeral #1@}
\makeatother

\begin{document}

\newcommand{\Cyc}{{\rm{Cyc}}}\newcommand{\diam}{{\rm{diam}}}
\newcommand{\Cays}{{\rm CayS}}
\newtheorem{thm}{Theorem}[section]
\newtheorem{pro}[thm]{Proposition}
\newtheorem{lem}[thm]{Lemma}
\newtheorem{exa}[thm]{Example}
\newtheorem{fac}[thm]{Fact}
\newtheorem{cor}[thm]{Corollary}
\theoremstyle{definition}
\newtheorem{ex}[thm]{Example}
\newtheorem{ob}[thm]{Observtion}
\newtheorem{remark}[thm]{Remark}
\newcommand{\bth}{\begin{thm}}
\renewcommand{\eth}{\end{thm}}
\newcommand{\bex}{\begin{ex}}
\newcommand{\eex}{\end{ex}}
\newcommand{\bre}{\begin{remark}}
\newcommand{\ere}{\end{remark}}

\newcommand{\bal}{\begin{aligned}}
\newcommand{\eal}{\end{aligned}}
\newcommand{\beq}{\begin{equation}}
\newcommand{\eeq}{\end{equation}}
\newcommand{\ben}{\begin{equation*}}
\newcommand{\een}{\end{equation*}}

\newcommand{\bpf}{\begin{proof}}
\newcommand{\epf}{\end{proof}}
\renewcommand{\thefootnote}{}
\newcommand{\sdim}{{\rm sdim}}

\def\beql#1{\begin{equation}\label{#1}}
\title{\Large\bf Perfect codes in Cayley sum graphs}
\author{{Xuanlong Ma$^{1}$, Kaishun Wang$^2$, Yuefeng Yang$^{3,}$\footnote{Corresponding author.}
}\\[15pt]
{\small\em $^1$School of Science, Xi'an Shiyou University, Xi'an 710065, China}\\
{\small\em $^2$Laboratory of Mathematics and Complex Systems {\rm (}Ministry of Education{\rm )},} \\
{\small\em School of Mathematical Sciences, Beijing Normal University, Beijing 100875, China
}\\
{\small\em $^3$School of Science, China University of Geosciences, Beijing 100083, China
}\\
}

 \date{}

\maketitle

\begin{abstract}
A subset $C$ of the vertex set of a graph $\Gamma$ is called a perfect code of $\Gamma$ if every vertex of $\Gamma$ is at distance no more than one to exactly one vertex in $C$. Let $A$ be a finite abelian group and $T$ a square-free subset of $A$. The Cayley sum graph of $A$ with respect to the connection set $T$ is a simple graph with $A$ as its vertex set, and two vertices $x$ and $y$ are adjacent whenever $x+y\in T$. A subgroup of $A$ is said to be a subgroup perfect code of $A$ if the subgroup is a perfect code of some Cayley sum graph of $A$. In this paper, we give some necessary and sufficient conditions for a subset of $A$ to be a perfect code of a given Cayley sum graph of $A$. We also characterize all subgroup perfect codes of $A$.
\end{abstract}


\medskip

\footnote{E-mail addresses: xuanlma@xsyu.edu.cn (X. Ma),
wangks@bnu.edu.cn (K. Wang), yangyf@cugb.edu.cn (Y. Yang).}
\section{Introduction}

In this paper, every group considered is finite, and
every graph considered is finite and simple.
For a graph $\Gamma$ with vertex set $V$,
a subset $C$ of $V$  is called a {\em perfect code} \cite{Kr}
of $\Gamma$ if every vertex of $\Gamma$ is at distance no more than one to exactly one vertex in $C$.
In other words, $C$ is a {\em perfect code} in $\Gamma$
provided that $C$ is independent in $\Gamma$ and
every vertex of $V\setminus C$ is adjacent to precisely one vertex of $C$.
In some references, a perfect code is also
called an {\em efficient dominating set} \cite{DYP,DYP2}
or {\em independent perfect dominating set} \cite{Le}.
Since the beginning of coding theory in the late 1940s, perfect codes have been one of the most important objects of study in
information theory. See, for example,
the two surveys \cite{He,Van} on perfect
codes and related definitions in the classical setting. Since the seminal paper of Biggs \cite{Big} and the fundamental work of Delsarte \cite{DD}, perfect codes in distance-transitive graphs and, in general, in distance-regular graphs and association
schemes have received
considerable attention in the literature.
Beginning with \cite{Kr}, a great
amount of work on perfect codes in general graphs has been produced. See, for example, \cite{DeS,Mo,Ze}.

For the past few years, perfect codes in Cayley graphs have attracted considerable attention, see, for example,
\cite{FHZ,Zh,Z15}. In \cite{HXZ18}, Huang, Xia and Zhou first introduced the concept of a perfect code of a group $G$.
A subset $C$ of $G$ is said to be a {\em perfect code} of $G$ if $C$ is a perfect code of some Cayley graph of $G$. In particular, a subgroup is said to be a {\em subgroup perfect code} of $G$ if the subgroup is also a perfect code of $G$. Also in \cite{HXZ18}, they gave a necessary and sufficient condition for a normal subgroup of a group $G$ to be a subgroup perfect code of $G$, and determined all the subgroup perfect codes of dihedral groups and some abelian groups. For more results on subgroup perfect codes of Cayley graphs, see
\cite{CWZ,MWWZ,ZZ20}.

In 1989, Chung~\cite{C89} first introduced the concept of a Cayley sum graph of an abelian group.
Let $A$ be an abelian group.
Given a subset $T$ of $A$, the {\em Cayley sum graph} (also called {\em addition Cayley graph}) of $A$ with respect to the connection set $T$,
denoted by ${\rm CayS}(A,T)$,
is a graph with $A$ as its vertex set, and two elements $x$ and $y$ are adjacent whenever $x+y\in T$.
An element $x$ of $A$ is said to be a {\em square} if $x = 2y$ for some $y\in A$. A subset of $A$ without squares is called a {\em square-free subset} of $A$.
Since every graph considered in this paper is simple, we always consider a simple Cayley sum graph ${\rm CayS}(A, T)$ of an abelian group $A$, that is,
the connection set $T$ should be square-free.
More explicitly, for a square-free subset $T$ of $A$, the {\em Cayley sum graph} ${\rm CayS}(A,T)$ of $A$ with respect to the connection set $T$ is a simple graph with $A$ as its vertex set, and two  elements $x$ and $y$ are adjacent whenever $x+y\in S$.
One can easily verify that ${\rm CayS}(A,T)$ is $|T|$-regular.
In \cite{Gry}, Grynkiewicz, Levb and Serra pointed out that, as the twins of the usual Cayley graphs, Cayley sum graphs are rather difficult to study, so
they received much less attention in the literature. For most results about Cayley sum graphs, see \cite{A07,Lev,Kon,DGM09,MW,CGW03}.

In \cite{MFW20}, the first author, Feng and the second author
studied the perfect codes of ${\rm CayS}(A, T)$, and
defined a subgroup perfect code of an abelian group by using Cayley sum graphs instead of Cayley graphs. More precisely, a subgroup of $A$ is said to be a {\em subgroup perfect code} of $A$ if the subgroup is a perfect code of some Cayley sum graph of $A$. Also, in \cite{MFW20}, the authors reduced
the problem of determining when a given subgroup of an abelian group is a perfect code to the case of abelian $2$-groups, and classified the abelian groups whose
all non-trivial subgroups are perfect codes.

In this paper, we study the perfect codes of a simple Cayley sum graph of an abelian group. For any abelian group $D$ of odd order, since $\{2g:g\in D\}=D$, it follows that every element of
$D$ is a square. Therefore, every simple Cayley sum graph of an abelian group $D$ of odd order is the empty Cayley sum graph ${\rm CayS}(D, \emptyset)$.
In order to study the non-trivial Cayley sum graphs of an abelian group, we always assume that the abelian group has even order.

Throughout the paper, $A$ denotes a finite abelian group of even order with identity element $0$.
The remainder of this paper is organized as follows.
In Section~\ref{2Sect}, we give some necessary and sufficient conditions for a subset of $A$ to be a perfect code of a given Cayley sum graph of $A$ (see Theorem~\ref{NN-thm1} and Corollary~\ref{NN-cor1}).
In Section~\ref{3Sect}, we determine the structure of a subgroup perfect code of $A$ (see Theorem~\ref{Ab-maimthm}),
which improves \cite[Theorem 3.1]{MFW20}, and we also give   some applications of Theorem~\ref{Ab-maimthm} (see Propositions~\ref{NNP1}, \ref{NNP2} and \ref{NNP3}).

\section{Perfect codes}\label{2Sect}

In this section, we study the perfect codes of a Cayley sum graph of $A$ and
give some necessary and sufficient conditions for a subset of $A$ to be a perfect code of a given Cayley sum graph of $A$ (see Theorem~\ref{NN-thm1} and Corollary~\ref{NN-cor1}).

For two subsets $B$ and $C$ of $A$, we write
$$
B\pm C=\{b\pm c: b\in B,~c\in C\},
$$
which is abbreviated by $b\pm C$ in the case where $B=\{b\}$.

\begin{lem}\label{NN-lemma1}
Take a subset $X$ of vertices in a Cayley sum graph $\Cays(A,T)$. Then every element of $A\setminus X$ is adjacent to at least one vertex of $X$ in ${\rm CayS}(A,T)$ if and only if
$$
A\setminus X\subseteq \bigcup_{t\in T}(t-X).
$$
\end{lem}
\proof Observe that $\bigcup_{t\in T}(t-X)$ consists of all vertices that are adjacent to some vertices of $X$. Hence, the desired result follows.
\qed

\begin{lem}\label{NN-lemma2}
Take a subset $X$ of vertices in a Cayley sum graph $\Cays(A,T)$. Then the following are equivalent:
\begin{itemize}
\item[{\rm (i)}] Every element of $A$ is adjacent to at most one element of $X$ in $\Cays(A,T)$;

\item[{\rm (ii)}] For each two distinct elements $t_1$ and $t_2$ in $T$, we have $(t_1-X)\cap (t_2-X)=\emptyset$;

\item[{\rm (iii)}] $(X-X)\cap(T-T)=\{0\}$.
\end{itemize}
\end{lem}
\proof Suppose that (i) holds. Assume that $(t_1-X)\cap (t_2-X)\neq\emptyset$ for some elements $t_1$ and $t_2$ in $T$. Then $t_1-x_1=t_2-x_2$ for some elements $x_1$ and $x_2$ in $X$. The fact that $T$ is square-free indicates that  $t_1-x_1\ne x_1$ and $t_1-x_1=t_2-x_2\ne x_2$. It follows that $t_1-x_1$ is adjacent to both $x_1$ and $x_2$,  which implies that $x_1=x_2$, and so $t_1=t_2$. Thus, (ii) is valid.

It is clear that (ii) implies (iii). Now suppose that (iii) holds. We only need to prove (i). Assume for a contrary that there exists an element $y$ in $A$ which is  adjacent to distinct elements $y_1$ and $y_2$ of $X$. Then $y+y_1\in T$ and $y+y_2\in T$, and so $y_1-y_2\in (X-X)\cap(T-T)$, a contradiction. This completes the proof.
\qed

\begin{lem}\label{NN-lemma3}
Take a subset $X$ of vertices in a Cayley sum graph $\Cays(A,T)$. Then the following are equivalent:
\begin{itemize}
\item[{\rm (i)}] $X$  is an independent set of $\Cays(A,T)$;

\item[{\rm (ii)}] For each $t\in T$, we have $X\cap (t-X)=\emptyset$;

\item[{\rm (iii)}] $(X+X)\cap T=\emptyset$.
\end{itemize}
\end{lem}
\proof
Now assume that (i) holds. Suppose to the contrary that
there exists $t\in T$ such that $X\cap (t-X)\ne\emptyset$.
Then there exist $x_1,x_2\in X$ such that $x_1=t-x_2$. Since $t$ is not a square, we have $x_1\ne x_2$.  It follows that $x_1$ is adjacent to $x_2$, this contradicts that  $X$ is an independent set. Thus, (ii) is valid.

Also, it is easy to see that (ii) implies (iii), and (iii) implies (i), the proof is complete.
\qed

\begin{thm}\label{NN-thm1}
For a Cayley sum graph $\Cays(A,T)$, write $T=\{t_1,\ldots,t_s\}$ and take a subset $X$ of $A$. The following are equivalent:
\begin{itemize}
\item[{\rm (i)}] $X$ is a perfect code of $\Cays(A,T)$;

\item[{\rm (ii)}] $\{X,t_1-X,\ldots,t_s-X\}$ is a partition of $A$;

\item[{\rm (iii)}] $|A|=|X|(s+1)$, $(X+X)\cap T=\emptyset$ and $(X-X)\cap(T-T)=\{0\}$.
\end{itemize}
\end{thm}
\proof
We first prove (i) implies (iii).
Suppose that (i) holds. By the definition of a perfect code, we have $|A|=|X|(s+1)$. Also, in view of Lemmas~\ref{NN-lemma2} and \ref{NN-lemma3}, we conclude  that (iii) holds.
Now by Lemmas~\ref{NN-lemma2} and \ref{NN-lemma3} again, it follows that (iii) implies (ii).

Now suppose that (ii) holds. It suffices to prove (i).
Note that
\begin{equation}\label{NN-ds1}
A=X\dot{\cup}(t_1-X)\dot{\cup}\cdots\dot{\cup}(t_s-X).
\end{equation}
Then $A\setminus X=\bigcup_{t\in T}(t-X)$.
By Lemma~\ref{NN-lemma1}, we have that every element of $A\setminus X$ is adjacent to at least one vertex of $X$ in $\Cays(A,T)$. Moreover, (\ref{NN-ds1}) implies that $(t_1-X)\cap (t_2-X)=\emptyset$
for each two distinct elements $t_1$ and $t_2$ in $T$.
It follows from Lemma~\ref{NN-lemma2} that every vertex in $A\setminus X$ is adjacent to exactly one vertex in $X$ in $\Cays(A,T)$. Notice that (\ref{NN-ds1}) also implies that $X\cap (t-X)=\emptyset$ for each $t\in T$. According to Lemma~\ref{NN-lemma3},  $X$ is an independent set of $\Cays(A,T)$, and so $X$ is a perfect code of $\Cays(A,T)$, as desired.
\qed

\bigskip

Two subsets $M$ and $N$ of $A$ are {\em supplementary}, denoted by $A=M\oplus N$, if each element $a$ of $A$ can be written as a unique manner $m+n$ with $m\in M$ and $n\in N$.

\begin{lem}{\rm (\cite[Proposition 2.1]{Vu})}\label{Vu1}
Let $M$ and $N$ be two subsets of $A$. Then $A=M\oplus N$ is equivalent to the conjunction of any two of the following conditions:
\begin{itemize}
\item[{\rm (a)}] $A=M+N$;

\item[{\rm (b)}] $(M-M)\cap(N-N)=\{0\}$;

\item[{\rm (c)}] $|A|=|M||N|$.
\end{itemize}
\end{lem}

\begin{cor}\label{NN-cor1}
For a Cayley sum graph $\Cays(A,T)$, take an inverse-closed subset $X$ in $A$. Then $X$ is a perfect code of $\Cays(A,T)$ if and only if $A=X\oplus T^0$, where $T^0=T\cup \{0\}$.
\end{cor}
\proof
Since $T$ is square-free, one obtains $0\notin T$, which implies that $|T^0|=|T|+1$.

We first prove the sufficiency. Suppose that $A=X\oplus T^0$.  Then by Lemma~\ref{Vu1}, we have $|A|=|X|(|T|+1)$ and  $(X-X)\cap(T^0-T^0)=\{0\}$.
Since $T-T\subseteq T^0-T^0$ and $T\subseteq T^0-T^0$, one gets $(X-X)\cap(T-T)=\{0\}$ and $(X-X)\cap T=\emptyset$. The fact that $X$ is inverse-closed implies that $(X+X)\cap T=\emptyset$. It follows from Theorem~\ref{NN-thm1} that $X$ is a perfect code of $\Cays(A,T)$, as desired.

We next prove the necessity. Suppose that $X$ is a perfect code of $\Cays(A,T)$. By Theorem~\ref{NN-thm1}, we have
$|A|=|X|(|T|+1)$, $(X+X)\cap T=\emptyset$ and $(X-X)\cap(T-T)=\{0\}$. Since $X$ is inverse-closed, one gets $(X-X)\cap T=\emptyset$, which implies $(-(X-X))\cap (-T)=\emptyset$. It follows that $(X-X)\cap (-T)=\emptyset$. Since
$T^0-T^0=T\cup (T-T) \cup (-T)$, we have that $(X-X)\cap(T^0-T^0)=\{0\}$. Now Lemma~\ref{Vu1} implies $A=X\oplus T^0$, as desired.
\qed

\bigskip

The cyclic group of order $n$ or the additive cyclic group of integers modulo $n$ is denoted by $\mathbb{Z}_n=\{0,1,\ldots,n-1\}$.
We use the following example to illustrate Corollary~\ref{NN-cor1}.

\begin{ex}
In $\mathbb Z_{12}$, let $X=\{0,3,6,9\}$ and $T=\{1,11\}$. Then $(X-X)\cap(T^0-T^0)=\{0\}$, and so $\mathbb Z_{12}=X\oplus T^0$.
Therefore, Corollary~\ref{NN-cor1} implies that
the Cayley sum graph $\Cays(\mathbb Z_{12},T)$ admits the perfect code $X$.
\end{ex}

If $\Cays(A,T)$ admits a perfect code, by Theorem~\ref{NN-thm1} we have that $|T|+1$ divides $|A|$.
Observe that $|A|$ must have divisor $2$, $|A|/2$ and $|A|$.
If $|T|+1=2$, since every $1$-regular graph has a perfect code, $\Cays(A,T)$ admits a perfect code.
If $|T|+1=|A|$, then $\Cays(A,T)$ is a complete graph, which admits a perfect code of size $1$ (in fact, in this case, $A$ is an elementary abelian $2$-group, see Corollary~\ref{yong-1}). Now we consider the case $|T|+1=|A|/2$.

\begin{cor}\label{NN-cor2}
A Cayley sum graph $\Cays(A,T)$ of valency $|A|/2-1$  admits a perfect code if $|T-T^0|<|A|$.
\end{cor}
\proof The condition $|T-T^0|<|A|$ indicates that there exists an element $a\in A$ such that $a\not\in T\cup(T-T)$. Let $X=\{0,a\}$. Then $X+X=\{0,a,2a\},X-X=\{0,a,-a\}$.
Since $T$ is square-free, we have $(X+X)\cap T=\emptyset$.
Since $a\not\in T-T$, one gets $-a\not\in T-T$. It follows that
$(X-X)\cap(T-T)=\{0\}$, and so the desired result follows from Theorem~\ref{NN-thm1}.
\qed

\bigskip

We conclude the section by the following example to illustrate Corollary~\ref{NN-cor2}.

\begin{ex}
Let $n$ be a positive even integer and let $T$ be a square-free subset of $\mathbb Z_n$. Then $T\subseteq \{1,3,\ldots,n-1\}$.
Suppose that $|T|=n/2-1$.
Let $\{a\}=\{1,3,\ldots,n-1\}\setminus T$. Then
$a\not\in T\cup(T-T)$. Thus, by Corollary~\ref{NN-cor2}, it follows that
a Cayley sum graph $\Cays(\mathbb Z_n,T)$ of valency $n/2-1$ has a perfect code. As a concrete example, the Cayley sum graph $\Cays(\mathbb Z_{12},\{1,3,5,7,9\})$ admits the perfect code $\{0,11\}$.
\end{ex}

\section{Subgroup perfect codes}
\label{3Sect}

In this section, we study the subgroup perfect codes of an abelian group of even order.
Our main result is Theorem~\ref{Ab-maimthm}
, which determines the structure of a subgroup perfect code of $A$  and improves \cite[Theorem 3.1]{MFW20}.
To state our main theorem, we prepare some basic notations.

Let $H$ be a subgroup of $A$.
The {\em index} of $H$ in $A$, denoted by $[A:H]$, is the number of right (or left) cosets of $H$ in $G$.
$H$ is said to be a {\em Hall $2'$-subgroup} of $A$ if $[A:H]$ is equal to the cardinality of a Sylow $2$-subgroup of $A$.
We use $A_2$ and $A_{2'}$ to denote the Sylow $2$-subgroup and Hall $2'$-subgroup of $A$, respectively.
Note that $A_2$ consists of the elements of $A$ each of whose order is a power of $2$, and $A_{2'}$ consists of the elements of $A$ with odd order. Particularly, if $A$ is a $2$-group, then $A=A_2$. As usual, we use $B \times C$ to denote the direct product (or direct sum) of two groups $B$ and $C$. So, we have $A=A_2 \times A_{2'}$.
By the fundamental theorem of finitely generated abelian groups, every finite abelian group is a direct product of some cyclic groups. Since $A=A_{2}\times A_{2'}$,
we may assume that
\begin{equation}\label{Ag-1}
A=\mathbb{Z}_{2^{m_1}}\times \mathbb{Z}_{2^{m_2}} \times \cdots
\times \mathbb{Z}_{2^{m_k}} \times A_{2'},
\end{equation}
where $m_i\ge 1$ for each $1\le i \le k$.
Observe that an element $(a_1,\ldots,a_k,a)$ in $A$ is not a square  if and only if there exists an odd integer in $\{a_1,\ldots,a_k\}$. In other words, an element $(a_1,\ldots,a_k,a)$ in $A$ is a square if and only if every element of $\{a_1,\ldots,a_k\}$ is even.
Remark that the whole group $A$ is a perfect code in the empty Cayley sum graph ${\rm CayS}(A, \emptyset)$.

\begin{thm}\label{Ab-maimthm}
Let $A$ be an abelian group as presented in {\rm (\ref{Ag-1})}, and let $H$ be a subgroup of $A$.
Then $H$ is a subgroup perfect code of $A$ if and only if
either $H$ is a subgroup isomorphic to
\begin{equation}\label{sz1}
\mathbb{Z}_{2^{m_1-1}}\times \mathbb{Z}_{2^{m_2-1}} \times \cdots
\times \mathbb{Z}_{2^{m_k-1}} \times A_{2'}
\end{equation}
or
$H$ has a non-square element.
\end{thm}

We use the following example to illustrate Theorem~\ref{Ab-maimthm}.

\begin{ex}
Let $A=\mathbb{Z}_2\times \mathbb{Z}_4\times \mathbb{Z}_3$. Then
$$H=\{(0,0,0),(0,2,0),(0,2,1),(0,2,2),(0,0,1),(0,0,2)\}$$
is the subgroup of $A$ which is isomorphic to $\mathbb{Z}_1\times \mathbb{Z}_2\times \mathbb{Z}_3$, and the set of all non-square elements of $A$ is
$$
A\setminus H=\{(1,x,y),(0,z,y): x\in \{0,1,2,3\},y\in \{0,1,2\},
z\in\{1,3\}\}.
$$
Therefore, Theorem~\ref{Ab-maimthm} implies that any subgroup perfect code of $A$ is either $H$ or a subgroup $K$ satisfying
$K\cap (A\setminus H)\ne \emptyset$. For example, the cyclic subgroup generated by any element belonging to $A\setminus H$ is a subgroup perfect code of $A$.
\end{ex}

In order to prove Theorem~\ref{Ab-maimthm}, we need some auxiliary results.
A {\em right transversal} (resp. {\em left transversal}) of a subgroup $H$ in $A$ is defined as a subset of $A$ which contains exactly one element in each right coset (resp. left coset) of $H$ in $A$.
In an abelian group, every right coset of any subgroup is also a left coset of the subgroup, for the sake of simplicity, we then use the term ``transversal'' to substitute for
``right transversal'' or ``left transversal''.
Our first result is the following proposition which gives some necessary and sufficient conditions for a subgroup of an abelian group to be a subgroup perfect code of the abelian group.
Remark that the following result has its own interest
and we only use one equivalent condition in the sequel.

\begin{pro}\label{Le-1}
Let $H$ be a subgroup of $A$.
The following are equivalent:
\begin{itemize}
\item[{\rm (i)}] $H$ is a subgroup perfect code of $A$;

\item[{\rm (ii)}] There exists a square-free subset $T\subseteq A$ such that $T^0$ is a transversal of $H$ in $A$;

\item[{\rm (iii)}] There exists a square-free subset $T\subseteq A$ such that
    $$[A:H]=|T|+1,~~H\cap (T\cup (T-T))=\{0\};$$

\item[{\rm (iv)}] There exists a square-free subset
$T\subseteq A$ such that $A=H\oplus T^0$;

\item[{\rm (v)}]  $H$ is a subgroup perfect code of any subgroup of $A$ which contains $H$;
\end{itemize}
\end{pro}
\proof
Combining Theorem~\ref{NN-thm1} and Corollary~\ref{NN-cor1}, one can verify that each two of (i), (ii), (iii) and (iv) are equivalent. We next prove that (i) and (v) are equivalent. It is clear that (v) implies (i).
Now suppose that $H$ is a subgroup perfect code of $A$.
Let $K$ be an arbitrary subgroup of $A$ with $H\subseteq K$.
It suffices to prove that $H$ is a subgroup perfect code of $K$.
Since (i) and (ii) are equivalent, there exists a square-free subset $T\subseteq A$ such that $T\cup\{0\}$ is a   transversal of $H$ in $A$. Let $T'=T^0\cap K$ where $T^0=T\cup \{0\}$. Then $0\in T'$ and $T'\setminus \{0\}$ is square-free.
Moreover, by the definition of a perfect code, one has
$$
K=A\cap K=(H+T^0)\cap K=H+(T^0\cap K)=H+T'.
$$
It follows that $T'$ is a transversal of $H$ in $K$. Since (i) and (ii) are equivalent again, $H$ is a subgroup perfect code of $K$.
\qed


\begin{lem}\label{Le-m1}
Let $H$ be a subgroup of $A$. If $H$ has a non-square element, then $H$ is a subgroup perfect code of $A$.
\end{lem}
\proof
Let $A$ have the form as presented in (\ref{Ag-1}).
Suppose that $(a_1,\ldots,a_k,a)\in H$ is a non-square element in $A$. Then there exists an odd integer in $\{a_1,\ldots,a_k\}$. In fact, without loss of generality,  it follows from (\ref{Ag-1}) that we may assume that $a_1$ is odd. Now let
$\{\alpha_0,\alpha_1,\ldots,\alpha_s\}$
be a transversal of $H$ in $A$, where $\alpha_0$ is the identity element of $A$ and $\alpha_i=(x_{i1},\ldots,x_{ik},x_{i})$ for each $1\le i \le s$.

For each $i\in\{1,\ldots,s\}$, if $\alpha_i$ is not a square,  then let $\beta_i=\alpha_i$; if $\alpha_i$ is a square, then every element of $\{x_{j1},\ldots,x_{jk}\}$ is even, and let
$$
\beta_i=(a_1,\ldots,a_k,a)+\alpha_i=(a_1+x_{i1},\ldots,a_k+x_{ik},a+x_i),
$$
which imply that $a_1+x_{i1}$ is odd in $\mathbb{Z}_{2^{m_1}}$ since $a_1$ is odd. It follows that $\beta_i$ is  a non-square element. Also, note that $\beta_i\in H+\alpha_i$. We conclude that
$\{\alpha_0,\beta_1,\ldots,\beta_s\}$
is a  transversal of $H$ in $A$, and $\{\beta_1,\ldots,\beta_s\}$ is square-free in $A$.
It follows from Proposition~\ref{Le-1} (ii) that $H$ is a subgroup perfect code of $A$.
\qed

\bigskip

Denote
$
A^{{\rm \Rmnum{2}}} = \{2a: a \in A\}.
$
Namely, $A^{{\rm \Rmnum{2}}}$ is the set of all squares of $A$.

\begin{lem}\label{Le-m2}
Let $H$ be a subgroup of $A$. Suppose that every element of $H$ is a square. Then $H$ is a subgroup perfect code of $A$ if and only if $H=A^{{\rm \Rmnum{2}}}$.
\end{lem}
\proof
Let $A$ have the form as presented in (\ref{Ag-1}).
We first prove the necessity. Suppose that $H$ is a subgroup perfect code of $A$. Note that $H\subseteq A^{{\rm \Rmnum{2}}}$. Assume the contrary, namely, there exists an element $(a_1,\ldots,a_k,a)\in A^{{\rm \Rmnum{2}}}\setminus H$. Then
every element of $\{a_1,\ldots,a_k\}$ is even. Also, by Proposition \ref{Le-1} (ii), we have that $H+(a_1,\ldots,a_k,a)$ contains at least one non-square element of $A$,
a contradiction since the sum of two squares is a square.

We now prove the sufficiency. Suppose that $H=A^{{\rm \Rmnum{2}}}$. Let
$\{\alpha_0,\alpha_1,\ldots,\alpha_s\}$
be a transversal of $H$ in $A$, where $\alpha_0$ is the identity element of $A$. It follows that $\{\alpha_1,\ldots,\alpha_s\}$ is square-free.
Now Proposition~\ref{Le-1} (ii) implies that $H$ is a subgroup perfect code of $A$.
\qed

\bigskip

We are now ready to prove Theorem~\ref{Ab-maimthm}.

\medskip

\noindent {\em Proof of Theorem~{\rm\ref{Ab-maimthm}}.}
We first claim that $H=A^{{\rm \Rmnum{2}}}$ if and only if
$H$ is a subgroup isomorphic to $(\ref{sz1})$.
In fact,
\begin{align*}
  A^{{\rm \Rmnum{2}}}& = \{2(a_1,\ldots,a_k,a):(a_1,\ldots,a_k,a)\in A\}  \\
   & = \{(2a_1,\ldots,2a_k,2a):(a_1,\ldots,a_k,a)\in A\}\\
   & = \{(b_1,\ldots,b_k,b):\text{$b_i$ is even and belongs to $\mathbb{Z}_{2^{m_i}}$ for all $1\le i \le k$, $b\in A_{2'}$}\}\\
   & \cong \mathbb{Z}_{2^{m_1-1}}\times \mathbb{Z}_{2^{m_2-1}} \times \cdots
\times \mathbb{Z}_{2^{m_k-1}} \times A_{2'}.
\end{align*}
Thus, the claim is valid. Now combining Lemmas~\ref{Le-m1} and \ref{Le-m2}, we complete the proof.
\qed

\bigskip

We next give a proof of \cite[Corollary 3.2]{MFW20} according to Theorem~\ref{Ab-maimthm}.


\begin{cor}{\rm (\cite[Corollary 3.2]{MFW20})}\label{Prop1}
An abelian group has a subgroup perfect code of odd order if and only if the Sylow $2$-subgroup of the group is an elementary abelian $2$-group and the subgroup perfect code is the Hall $2'$-subgroup of the group.
\end{cor}
\proof
Let $A$ be an abelian group as presented in {\rm (\ref{Ag-1})}.
The proof of the sufficiency is straightforward by Theorem~\ref{Ab-maimthm}.
We now prove the necessity.
Suppose that $H$ is a subgroup perfect code of $A$ and has odd order. Then $H$ has no non-square elements.
Theorem~\ref{Ab-maimthm} implies that $H$ is a subgroup isomorphic to $(\ref{sz1})$. It follows that $m_i=1$ for all $1\le i \le k$, and so the Sylow $2$-subgroup of $A$ is  elementary abelian and $H$ is the Hall $2'$-subgroup of $A$, as desired.
\qed

\bigskip

The following is immediate by Corollary~\ref{Prop1}.

\begin{cor}\label{yong-1}
The subgroup consisting of the identity element of an abelian group
is a subgroup perfect code if and only if the abelian group is an elementary abelian $2$-group.
\end{cor}

We next give a necessary condition for a subgroup of an abelian group to be a subgroup perfect code of the abelian group.

\begin{pro}\label{FPo1}
Let $H=H_2\times H_{2'}$ be a subgroup of $A=A_2\times A_{2'}$.
If $H$ is a subgroup perfect code of $A$, then $H_2$ is a subgroup perfect code of $A_2$.
\end{pro}
\proof If $H=A$, then $H_2=A_2$, and so $H_2$ is a subgroup perfect code of $A_2$ from Theorem \ref{Ab-maimthm}, as desired. Thus, in the following, we may assume that $H\ne A$. In view of Proposition \ref{Le-1} (ii), the group $A_2\times A_{2'}$ has a square-free subset $T:=\{(x_1,y_1),(x_2,y_2),\ldots,(x_s,y_s)\}$ such that $T\cup\{(0,0)\}$ is a transversal of $H_2\times H_{2'}$ in $A_2\times A_{2'}$. Note that $H_2\subseteq A_2$ and $H_{2'}\subseteq A_{2'}$.
It follows that $\{x_1,x_2,\ldots,x_s\}$ is a square-free subset of $A_2$. Notice that $\bigcup_{i=0}^s(H_2+x_i)=A_2$, where $x_0=0$. So, we may assume that $T':=\{x_1',x_2',\ldots,x_{s'}'\}$ is a subset of $\{x_1,x_2,\ldots,x_s\}$ such that $T'$ is square-free and $T'\cup\{0\}$ is a transversal of $H_2$ in $A_2$.
Note that (i) and (ii) in Proposition \ref{Le-1} are equivalent. It follows that $H_2$ is a subgroup perfect code of $A_2$, as desired.
\qed

\medskip

By Corollary~\ref{Prop1}, it is easy to see that the converse of
Corollary~\ref{FPo1} is not true.
Finally, as some applications of Theorem \ref{Ab-maimthm}, we determine the subgroup perfect codes of three families of abelian groups. The first result is obtained by applying
Theorem~\ref{Ab-maimthm} to a cyclic group, which determines all
subgroup perfect codes of a cyclic group of even order.

\begin{pro}{\rm (\cite[Theorem 3.7]{MFW20})}\label{NNP1}
Let $A=\mathbb{Z}_{2^n}\times \mathbb{Z}_{m}$ be a cyclic group of even order, where $n\ge 1$ and $m$ is an odd integer. Then a subgroup of $A$ is a subgroup perfect code if and only if
the subgroup is isomorphic to either $\mathbb{Z}_{2^{n-1}}\times \mathbb{Z}_{m}$ or $\mathbb{Z}_{2^{n}}\times \mathbb{Z}_{m'}$, where $m'\mid m$.
\end{pro}

\begin{pro}\label{NNP2}
Let $A=\mathbb{Z}_2^n\times A_{2'}$, where $n\ge 1$. Then a subgroup of $A$ is a subgroup perfect code if and only if either the subgroup has even order or the subgroup is $A_{2'}$.
\end{pro}
\proof
The necessity follows trivially from Corollary~\ref{Prop1}.
We now prove the sufficiency. By Corollary~\ref{Prop1}, we only need to prove that if a subgroup of $A$ has even order, then
the subgroup is a subgroup perfect code. Now suppose that $H$ is a subgroup of $A$ and has even order. Taking $a\in H$ with order $2$, we have $a\in A_2$. Since $A_2$ is an elementary abelian $2$-group, $A$ has no elements of order $4$, which implies that $a$ is a non-square element. Lemma~\ref{Le-m1} implies that $H$ is a subgroup perfect code, as desired.
\qed

\begin{pro}\label{NNP3}
Let $A=\mathbb{Z}_2^n\times \mathbb{Z}_4 \times A_{2'}$, where $n\ge 1$. Then a subgroup of $A$ is a subgroup perfect code if and only if it is not isomorphic to one of $E_1^n\times \mathbb{Z}_2 \times A_{2'}'$ and $A_{2'}''$, where $E_1$ is the group consisting of the identity element of $\mathbb{Z}_2$, $A_{2'}'$ is a proper subgroup of $A_{2'}$, and $A_{2'}''$ is a subgroup of $A_{2'}$.
\end{pro}
\proof
We first prove the necessity. Let $H$ be a subgroup perfect code of $A$.
Clearly, by Corollary~\ref{Prop1}, $H$ is not isomorphic to a subgroup of $A_{2'}$.
Suppose to the contrary that $H\cong E_1^n\times \mathbb{Z}_2 \times A_{2'}'$, where $E_1$ is the group consisting of the identity element of $\mathbb{Z}_2$ and $A_{2'}'$ is a proper subgroup of $A_{2'}$. Then any element of $H$ has the form
\begin{center}
$(\underbrace{0,\ldots,0}_n,0,a)$ or
$(\underbrace{0,\ldots,0}_n,2,a)$,
\end{center}
where $a\in A_{2'}'$. It follows that every element of $H$ is a square.
By Theorem~\ref{Ab-maimthm}, it follows that $A_{2'}'=A_{2'}$,
which contradicts that $A_{2'}'$ is a proper subgroup of $A_{2'}$.

We next prove the sufficiency. Suppose that a subgroup $H$ of $A$ is not isomorphic to one of $E_1^n\times \mathbb{Z}_2 \times A_{2'}'$ and $A_{2'}''$. If $H\cong E_1^n\times \mathbb{Z}_2 \times A_{2'}$, then Theorem~\ref{Ab-maimthm} implies that $H$ is a subgroup perfect code, as desired. Now suppose that $H\ncong E_1^n\times \mathbb{Z}_2 \times A_{2'}$. Then by the structure of $A$, we conclude that $H$ has a non-square element, and so
$H$ is a subgroup perfect code by Lemma~\ref{Le-m1}, as required.
\qed

\bigskip
\noindent \textbf{Acknowledgements}~~
X. Ma was supported by the
National Natural Science Foundation of China
(11801441, 61976244), the Natural Science Basic Research Program of Shaanxi (Program No. 2020JQ-761)
and the Young Talent fund of University Association for Science and Technology in Shaanxi, China (20190507). K. Wang is supported by the National Natural Science
Foundation of China (11671043). Y. Yang was supported by the Fundamental Research Funds for the Central Universities
(2652019319).

\end{document}